\documentclass[proc]{edpsmath}
\usepackage{amsthm}
\usepackage{enumitem}
\theoremstyle{remark}
\usepackage{graphicx}
\usepackage{float}

\begin{document}

\title{Stability theory for some scalar finite difference schemes : Validity of the modified equations approach}
\author{Firas Dhaouadi}\address{Universit\'{e} Paul Sabatier, Institut de Math\'{e}matiques de Toulouse. e-mail :  dhaouadi@insa-toulouse.fr}
\author{Emilie Duval}\address{Universit\'{e} Grenoble Alpes, Laboratoire Jean Kuntzmann. e-mail: emilie.duval@inria.fr}
\author{Sergey Tkachenko}\address{Aix-Marseille Universit\'{e}, CNRS, IUSTI, UMR 7343. e-mail: sergey.tkachenko@univ-amu.fr}
\author{Jean-Paul Vila}\address{Institut de Math\'{e}matiques de Toulouse, INSA Toulouse. e-mail: vila@insa-toulouse.fr}
%
%
\begin{abstract} In this paper, we discuss some limitations of the modified equations approach as a tool for stability analysis for a class of explicit linear schemes to scalar partial derivative equations. We show that the infinite series obtained by Fourier transform of the modified equation is not always convergent and that in the case of divergence, it becomes unrelated to the scheme. Based on these results, we explain when the stability analysis of a given truncation of a modified equation may yield a reasonable estimation of a stability condition for the associated scheme. We illustrate our analysis by some examples of schemes namely for the heat equation and the transport equation. \end{abstract}
\begin{resume}
Dans ce papier, nous examinons quelques limites des \'{e}quations modifi\'{e}es, comme outil pour l'analyse de la stabilit\'{e} d'une certaine classe de sch\'{e}mas num\'{e}riques explicites lin\'{e}aires pour des \'{e}quations aux d\'{e}riv\'{e}es partielles scalaires. Nous montrons que la s\'{e}rie obtenue par transform\'{e}e de Fourier de l'\'{e}quation modifi\'{e}e n'est pas toujours convergente et que dans ce cas, son comportement n'est plus li\'{e} \`{a} celui du sch\'{e}ma. A partir de ces r\'{e}sultas, nous expliquons quand l'analyse de stabilit\'{e} d'une troncature donn\'{e}e peut donner des conditions de stabilit\'{e} raisonnables. Cette analyse est illustr\'{e}e par quelques exemples, notamment des sch\'{e}mas pour l'\'{e}quation de la chaleur  et pour l'\'{e}quation de transport.
 \end{resume}
\maketitle

\section*{Introduction}

Modified equations have been the subject of investigations and debates in numerical analysis. And yet, despite being relatively easy to establish, their use has been regarded with a lot of skepticism due to the lack of theoretical justification and the trickiness of their analysis. The first use of this technique for stability purposes is due to Hirt \cite{hirt1968heuristic}. He provided the first practical examples for which modified equations give relevant information, in a heuristic manner. The pionneering work of Warming and Hyett \cite{Fluid1974} introduced how to obtain these equations in the general case and mostly clarified the link between the stability of the scheme and the modified equations. They provided a simple and efficient way to obtain the exact Von Neumann stability conditions from a finite number of the modified equation coefficients, for linear scalar schemes. An alternate way to obtain modified equations was introduced in \cite{carpentier1997derivation}. The method permits to obtain the same modified equation through a series expansion of a an explicitly known function rather than the elimination technique of \cite{Fluid1974}. Other works that tried to tackle stability analysis through modified equations include \cite{li2013neumann,tyler1971heuristic}.

In this paper, we investigate the modified equations for some finite difference schemes, in an attempt to make it clear why this technique  often fails to provide relevant information on stability. First, we present a reminder on how to obtain modified equations for a given linear scalar scheme. We explain the basics of the heuristic stability theory and its limitations through some examples. In the second part, we present the technical framework and tools needed for stability analysis. We then clarify the mathematical reasons behind the frequent failures of the heuristic stability theory by investigating the convergence of the Fourier transform of the modified equation. We show that the latter is conditionally convergent and only then does it give significant results. Then, we compare the stability conditions of truncated modified equations with the corresponding scheme in the region of convergence. Finally, we present some examples to justify our analysis.
\section{On modified equations and heuristic stability}
\subsection{Obtaining the equations}
Consider as an example, the linear scalar transport equation given by : 
\begin{equation}
    \frac{\partial \tilde{u}}{\partial t} + c \frac{\partial \tilde{u}}{\partial x} = 0
    \label{eq:trans}
\end{equation}
with positive velocity $c$. In order to solve this equation with a finite difference scheme, we first introduce a uniform grid of points defined as usual by $(x_j=j\Delta x,t^n=n \Delta t)$ and we denote by $u_j^n=u(x_j,t^n)$ the value of the numerical solution in the corresponding grid point. Under these notations, take for instance the upwind Euler scheme for equation (\ref{eq:trans}) : 
\begin{equation}
    \frac{u_j^{n+1}-u_j^n}{\Delta t} +c \frac{u_j^n-u_{j-1}^n}{\Delta x}=0
    \label{eq:scheme_trans}
\end{equation}
Generally, in order to get relevant information on the consistency of the scheme \cite{richtmyer1994difference} or its convergence rate, we assume the existence of a smooth, infinitely differentiable numerical solution $u(x,t)$ that satisfies  $u(x_j,t^n)=u_j^n$ in each grid point. Provided this solution, we expand each term of the scheme  in Taylor series in the vicinity of $(x_j,t^n)$ which leads to an equation containing an infinite number of partial derivatives : 
\begin{equation}
    \frac{\partial u}{\partial t} + \frac{\Delta t}{2} \frac{\partial^2 u}{\partial t^2} + \frac{\Delta t^2}{6} \frac{\partial^3 u}{\partial t^3} + \ldots = -c \frac{\partial u}{\partial x} + \frac{c \Delta x}{2} \frac{\partial^2 u}{\partial x^2} + \frac{c \Delta x^2}{6} \frac{\partial^3 u}{\partial t^3} + \ldots
    \label{eq:mod1}
\end{equation}
This equation as is, is sufficient to prove consistency. In fact, we can clearly see that the in the limit $\Delta t\rightarrow 0$ and $\Delta x\rightarrow 0$ we recover the original transport equation. However, for further analysis of the numerical effects induced by the scheme, it would be more intuitive to consider an evolution equation with only space derivatives. To do that, we use the elimination procedure introduced by Warming and Hyett in \cite{Fluid1974}, that is we repeatedly use linear combinations of the equation (\ref{eq:mod1}) and its derivatives in order to eliminate higher order time derivatives and finally obtain the equation : 
\begin{equation}
    \frac{\partial u}{\partial t} + c \frac{\partial u}{\partial x} = c \frac{\Delta x}{2} \left(1- c \frac{\Delta t}{\Delta x} \right) \frac{\partial^2u}{\partial x^2} -\frac{\Delta x^2}{6}\left(1- c \frac{\Delta t}{\Delta x} \right)\left(1- 2c \frac{\Delta t}{\Delta x} \right)\frac{\partial^3u}{\partial x^3} \ldots
\end{equation}
This is called the modified equation associated with the upwind Euler scheme for the transport equation (\ref{eq:trans}). It is worth noting that this is not a partial derivative equation in the conventional sense as it does not have an order or a finite amount of partial derivatives. For a solution to exist, one also needs a  proper definition of infinitely many boundary conditions. Therefore, we restrict our analysis to solutions that are periodic \cite{Fluid1974}.

\subsection{Heuristic stability and limitations}
Besides proving consistency, the modified equation  quantifies explicitly the additional numerical effects. For example, it tells that up to first order in $\Delta x$ and $\Delta t$, the numerical solution rather satisfies a convection diffusion equation with a numerical dissipation coefficient given by $c \frac{\Delta x}{2} \left(1- c \frac{\Delta t}{\Delta x} \right)$. Note however that the sign of the latter can be negative in which case this equation becomes unstable and admits solutions that grow exponentially in time rather than decay. It follows from this that a necessary and sufficient condition of stability for the solutions of this convection-diffusion equivalent is : 
\begin{equation}
    c\frac{\Delta t}{\Delta x}\leq 1
\end{equation}
which turns out to be exactly the CFL condition for the scheme (\ref{eq:scheme_trans}). This gives a lot of potential for the modified equations to be a practical tool for stability analysis. Although the analysis is completely heuristic and has no rigorous foundation, its results make sense for a large class of schemes. However, there are also many examples for which this analysis fails to provide any practical stability condition. Consider for example the one dimensional heat equation : 
\begin{equation}
    \frac{\partial u}{\partial t}-\alpha \frac{\partial^2 u}{\partial x^2} = 0
\end{equation}
where $\alpha>0$ is the diffusion coefficient. We discretize this equation using centered finite differences : 
\begin{equation}
    \frac{u_i^{n+1} - u_i^n}{\Delta t} -\alpha \frac{u_{i+1}^n-2u_i^n+u_{i-1}^n}{\Delta x^2}= 0
    \label{eq:scheme_heat}
\end{equation}
Using the same elimination procedure we obtain the following modified equation up to $4th$ order : 
\begin{equation}
    \frac{\partial u}{\partial t}-\alpha \frac{\partial^2u}{\partial x^2} =  \frac{\alpha \Delta x^2}{12}\left(1-6\frac{\Delta t}{\Delta x^2} \right) \frac{\partial^4u}{\partial x^4} + \ldots
    \label{eq:mod_heat}
\end{equation}
If we choose to stop at the first non-zero truncation term as previously, the conclusions are already less straightforward as we have a non-vanishing second order term that comes from the heat equation itself and a $4th$ order term that comes from the numerical effects. Therefore in order to look into stability of this equation, let us shift to Fourier space. Let $v(k,t)$
be the Fourier transform of $u(x,t)$ then under these notations, equation (\ref{eq:mod_heat}) yields : 
\begin{equation}
\frac{\partial v}{\partial t} = -\alpha k^2\left(1-\frac{k^2\Delta x^2}{12}\left(1-6\frac{\Delta t}{\Delta x^2}\right)\right)v 
\end{equation}
It is reasonable to only consider the wavenumbers that are bound by $\left\vert k\Delta x \right\vert\leq \pi$, since these are the only wavenumbers admissible by the discrete mesh . In this setting, it is easy to verify that : 
\begin{equation}
    -\alpha k^2\left(1-\frac{k^2\Delta x^2}{12}\left(1-6\frac{\Delta t}{\Delta x^2}\right)\right) \leq 0 \quad \forall k \in \left[-\frac{\pi}{\Delta x},\frac{\pi}{\Delta x} \right]
\end{equation}
which implies that the considered truncation is unconditionally stable for all admissible $\Delta t$ and $\Delta x$. This result is obviously erroneous as it is well known that the Von Neumann stability analysis proves that a necessary and sufficient stability condition for the scheme (\ref{eq:scheme_heat}) is :
\begin{equation}
  \alpha \frac{\Delta t}{\Delta x^2}\leq \frac{1}{2}
\end{equation}
Thus, the heuristic analysis of this truncation of the modified equation failed to provide any practical stability condition for the scheme and truncating the equation at higher orders does not seem to do any better. This is one of the examples that demonstrates limitations of the method.
In what follows, we first introduce all the necessary notations and setting for stability analysis before attempting to explain why the stability of the truncation is sometimes incoherent with the stability of the scheme.
\section{Theory of stability through modified equations}

\subsection{Notations and assumptions}
We will consider for our analysis partial derivative equations that are linear, first order in time and of arbitrary order in space : 
\begin{equation}
\frac{\partial \tilde u}{\partial t} + \sum_{p=1}^{P} A_{p}\frac{\partial^p \tilde u}{\partial x^p} = 0
\label{eq:pde}
\end{equation}
where $A_p$ are constants. We consider explicit in time schemes, that are consistent with equation (\ref{eq:pde}) and can be written as : 
\begin{equation}
u_j^{n+1} = u_j^n + \Delta t \sum_{p=-n_l}^{n_r} b_{p}\left(\Delta x\right) u_{j+p}^n
\label{eq:scheme}
\end{equation}
where $n_l$ and $n_r$ are the number of mesh points to the left and to the right of $x_j$ respectively, used in every iteration. The coefficients $b_p$ verify 
\begin{equation} \sum_{p=-n_l}^{n_r} b_{p}(\Delta x)=0
\label{eq:consis}
\end{equation}
for consistency purposes, so that constant solutions, which are solution of the pde (\ref{eq:pde}) remain as such for the scheme. Sometimes it is preferable to cast the scheme (\ref{eq:scheme}) into an equivalent form : 
\begin{equation}
u_j^{n+1} = u_j^n + \frac{\Delta t}{\Delta x^q} \sum_{p=-n_l}^{n_r} B_{p}\left(\Delta x\right) u_{j+p}^n
\label{eq:scheme_lamdap}
\end{equation}
where $q$ is the highest power of $1/\Delta x$ present in the summation. Very often, when using standard finite differences, the coefficients $b_p$ are polynomials in $1/\Delta x$ and $q$ is the highest degree among these polynomials which is frequently equal to the order $P$ of the pde (\ref{eq:pde}). This formulation can be practical for stability analysis since most stability conditions for explicit schemes are given by bounds on the quantity $\lambda_q = \Delta t/\Delta x^q$ in the limit $\Delta x\rightarrow 0$ and $\Delta t\rightarrow 0$. Therefore, setting $\lambda_q$ as a non-vanishing parameter is a reasonable assumption that permits to reduce the number of free small parameters, that is we can take $\Delta t=g_q(\Delta x)=\lambda_q \Delta x^q$ and consider instead $\lambda_q$ and $\Delta x$ as free independent parameters.
In what follows, in order to perform stability analysis in Fourier space, we consider a space continuous counterpart of the scheme (\ref{eq:scheme_lamdap}) : 
\begin{equation}
u^{n+1}(x) = u^n (x)+ \lambda_q \sum_{p=-n_l}^{n_r} B_{p}\left(\Delta x\right)u^n (x+p\Delta x)
\label{eq:scheme_cont}
\end{equation}
Let $v(k,t)$ be the Fourier transform in space of $u(x,t)$. If we take $\theta=k\Delta x$, then the Fourier transform of equation (\ref{eq:scheme_cont}) yields : 
\begin{equation}
    v^{n+1}(k) = \left( 1+\lambda_q \sum_{p=-n_l}^{n_r} B_{p}(\Delta x)\mathrm{e}^{\mathrm{i}p\theta}\right)v^n(k) = S(\theta,\lambda_q,\Delta x)v^n(k)
\end{equation}
Since we will be operating most of the time in Fourier space, it seems necessary to recall consistency of the scheme in the same setting. Therefore if we assume that the exact solution to the pde (\ref{eq:pde}) satisfies in Fourier space: 
\begin{equation}
\tilde{v}(k,t+\Delta t) = \exp\left(\Delta t\sum_{p=1}^{P}(\mathrm{i}k)^p A_{p}\right)\tilde{v}(k,t) =\tilde{G}(\theta,\Delta t) \tilde{v}(k,t) 
\end{equation}
then the scheme (\ref{eq:scheme}) is consistent with the pde (\ref{eq:pde}) to order $s$ if and only if \cite{dautray2012mathematical} :
\begin{equation}
\left\vert \frac{\tilde{G}(\theta,\Delta t)-S(\theta,\lambda_q,\Delta x)}{\Delta t} \right\vert = \mathcal{O}(\Delta t^s)
\end{equation}
Lastly, we denote the modified equation associated to the scheme (\ref{eq:scheme}) by : 
\begin{equation}
u_t = \sum_{p=1}^{\infty} \mu_p(\lambda_q,\Delta x) \frac{\partial^pu}{\partial x^p}
\label{eq:mod}
\end{equation}
where $\mu_p$ are constants depending on $\Delta x$ and $\lambda_q$. 
\subsection{Fourier stability analysis}
In this part, we focus on the link between the modified equation and the scheme \cite{Fluid1974}. The amplification factor of the scheme is none other than the modulus of $S(\theta,\lambda_q,\Delta x)$. Now, in order to recover an equivalent expression in the continuous time for the modified equation, we apply the Fourier transform to  (\ref{eq:mod}). This implies that $v(k,t)$ satisfies the differential equation: 
\begin{equation}
\frac{dv}{dt}= \left( \sum_{p=1}^{\infty} (ik)^p\mu_p(\lambda_q,\Delta x) \right)v(k,t) =  \left( \sum_{p=1}^{\infty} \alpha_p(\lambda_q, \Delta x) \theta^p \right)v(k,t) =G(\theta,\lambda_q,\Delta x)v(k,t)
\label{eq:Fourier_mod}
\end{equation} 
where $\alpha_p(\lambda_q,\Delta x) = \mu_p(\lambda_q,\Delta x)i^p/\Delta x^p$. For convenience, we will say that the modified equation is stable if, for an initial condition $v(k,t=0)=v_0(k)$, the solution to the Cauchy problem : 
\begin{equation}
\left\{\begin{array}{l}
\displaystyle\frac{dv}{dt} = G(\theta,\lambda_q,\Delta x)v(k,t) \\
v(k,0)=v_0(k)
\end{array} \right.
\label{eq:cauchy}
\end{equation}
that is given by:
\begin{equation}
v(k,t)=\mathrm{e}^{tG(\theta,\lambda_q,\Delta x)}v_0(k)
\end{equation}
remains bounded $\forall t <T$, where $T>0$.
Given this solution, one can obviously write : 
\begin{equation}
v(k,t+\Delta t) = \mathrm{e}^{\Delta t G(\theta,\lambda_q,\Delta x)}v(k,t)
\end{equation}
Now, since the solution to the scheme (\ref{eq:scheme}) with the same initial condition is also an exact solution to the modified equation (\ref{eq:mod}) \cite{Fluid1974}, uniqueness of this solution gives :
\begin{equation}
\mathrm{e}^{\Delta t G(\theta,\lambda_q,\Delta x)} = \mathrm{e}^{\lambda_q\Delta x^q G(\theta,\lambda_q,\Delta x)} =S(\theta,\lambda_q,\Delta x)
\label{eq:equality}
\end{equation}
It follows from that the proposition : 
\begin{prpstn}
\label{pr:1}
For any scheme that writes as (\ref{eq:scheme}) and that is consistent with (\ref{eq:pde}), there exists a positive number $\theta_m$ that depends only on $\lambda_q$ and $\Delta x$ such that, $\forall \theta \in \left]-\theta_m,\theta_m\right[$ the expansion 
\begin{equation*}
    G(\theta,\lambda_q,\Delta x) = -\sum_{p=1}^{\infty}\frac{(1-S(\theta,\lambda_q,\Delta x))^p}{p\lambda_q\Delta x^q} 
\end{equation*}
holds and the series $\sum_{p=1}^{\infty} \alpha_p(\lambda_q,\Delta x) \theta^p$ converges to $\frac{1}{\lambda_q \Delta x^q} \ln(S(\theta,\lambda_q,\Delta x))$.
\end{prpstn}

\begin{proof}
The scheme (\ref{eq:scheme}) is consistent with the pde (\ref{eq:pde}) and so we have $\left\vert S(0,\lambda_q,\Delta x)-1\right\vert=0 < 1$. Since $S$ is a continuous function with respect to $\theta$, $\exists \theta_m(\lambda_q,\Delta x)>0$ such that $\left\vert S(\theta,\lambda_q,\Delta x)-1\right\vert < 1\  \forall \theta \in ]-\theta_m,\theta_m[$, and consequently, the principal logarithm $\ln(S(\theta,\lambda_q,\Delta x))$ defined by the series expansion : 
\begin{equation}
\ln(S(\theta,\lambda_q,\Delta x)) = \ln(1-(1-S(\theta,\lambda_q,\Delta x)))=\sum_{p=1}^{\infty}-\frac{(1-S(\theta,\lambda_q,\Delta x))^p}{p}  
\end{equation}
is convergent and $\mathrm{e}^{\ln(S(\theta,\lambda_q,\Delta x))}=S(\theta,\lambda_q,\Delta x)$. This proves that the function:
\begin{equation}
G(\theta,\lambda_q,\Delta x) = \frac{1}{\lambda_q\Delta x^q} \ln(S(\theta,\lambda_q,\Delta x))
\end{equation}
is a solution of equation (\ref{eq:equality}) for $\theta\in]-\theta_m,\theta_m[$. The exponential function is locally invertible in the vicinity of $0$ and so this expansion is unique in this vicinity \cite{carpentier1997derivation} and therefore $\forall \theta \in ]-\theta_m,\theta_m[$ . Moreover, since $S(\theta,\lambda_q,\Delta x)$ is entire (as a finite sum of exponential functions) then further expanding $S$ into power series of $\theta$ for $\theta \in ]-\theta_m,\theta_m[$ finally yields: 
\begin{equation}
G(\theta,\lambda_q,\Delta x) =  \sum_{p=1}^{\infty} \alpha_p \theta^p  = \frac{1}{\lambda_q \Delta x^q} \ln(S(\theta,\lambda_q,\Delta x))
\label{eq:G=ln}
\end{equation}
and hence, $\sum_{p=1}^{\infty} \alpha_p(\lambda_q,\Delta x) \theta^p$ is the series expansion of $\frac{1}{\lambda_q \Delta x^q} \ln(S(\theta,\lambda_q,\Delta x))$. This concludes our proof. 
\end{proof}
\begin{rmrk}
Equality (\ref{eq:G=ln}) also implies that the coefficients $\alpha_p(\lambda_q,\Delta x)$ are given by : 
\begin{equation}
    \alpha_p(\lambda_q,\Delta x) = \left.\frac{\partial^p }{\partial \theta^p}\left(\frac{\ln(S(\theta,\lambda_q,\Delta x))}{p!\lambda_q\Delta x^q}\right)\right\vert_{\theta=0}
\end{equation}
\end{rmrk}
\begin{rmrk} $\theta_m$ is not the radius of convergence of the series $G(\theta,\lambda_q,\Delta x)$. If we denote by $R$ the radius of convergence, then we have $R\geq \theta_m$. The exact radius is not trivial to find in practice in the general case, since $G(\theta,\lambda_q,\Delta x)$ is a series expansion of a composite function. However it is sometimes possible to give estimates or exact values of the radius for some examples as will be shown later. 
\end{rmrk}
\begin{rmrk}
In practice, when looking for convergence, we are mainly searching for constraints in $\lambda_q$ and $\Delta x$, for which the series $G(\theta,\lambda_q,\Delta x)$ is convergent  $ \forall \theta \in [-\pi,\pi]$, that is we want $R>\pi$. A sufficient condition is $\theta_m>\pi$. 
\end{rmrk}

So far, we have shown that for fixed parameters $\Delta x$ and $\lambda_q$, the series $G(\theta,\lambda_q,\Delta x)$ converges if $\theta \in ]-\theta_m,\theta_m[$. This being a sufficient condition of convergence, we do not know what happens for $\left\vert \theta\right\vert\geq \theta_m$. On the other hand it is worth noting that, for $\left\vert\theta\right\vert > R$, the series is divergent and equality (\ref{eq:equality}) does not hold anymore. Therefore, the modified equation stability is not linked to that of the scheme for $\left\vert\theta\right\vert>R$.

\subsection{Scheme stability domain and series convergence domain}
Generally, in the Von Neumann setting, stability conditions of the scheme are given by constraints linking $\lambda_q$ and $\Delta x$ so that the inequality 
\begin{equation}
 \left\vert S(\theta,\lambda_q,\Delta x)\right\vert \leq 1 \quad \forall \theta \in \left[-\pi,\pi\right]
\end{equation}
is verified. These constraints define a region of stability $\mathcal{R}_s$ in the $(\lambda_q,\Delta x)$ plane, that is : 
\begin{equation}
    \mathcal{R}_s = \left\{ (\lambda_q,\Delta x) \in \mathbb{R}_+^2 :    \forall \theta \in \left[-\pi,\pi\right] :  \left\vert S(\theta,\lambda_q,\Delta x)\right\vert \leq 1 \right\} 
\end{equation}
In the same manner, we can define a region $\mathcal{R}_c$ of the same plane, in which the series $G(\theta,\lambda_q,\Delta x)$ converges :
\begin{equation}
\mathcal{R}_c = \left\{ (\lambda_q,\Delta x) \in \mathbb{R}_+^2 : \forall \theta \in \left[-\pi,\pi\right] : \left \vert G(\theta,\lambda_q,\Delta x) \right\vert < \infty \right\}  
\end{equation}
which is also equivalent to :
\begin{equation}
  \mathcal{R}_c  = \left\{ (\lambda_q,\Delta x) \in \mathbb{R}_+^2 : R(\lambda_q,\Delta x) \geq \pi\right\}
\end{equation}
In practice, it is not always possible to exactly determine $\mathcal{R}_c$. It is however easier to explicitly find a subset of this region, that is :
\begin{equation}
\Omega_c =   \left\{ (\lambda_q,\Delta x) \in \mathbb{R}_+^2 :    \forall \theta \in \left[-\pi,\pi\right] :  \left\vert 1- S(\theta,\lambda_q,\Delta x)\right\vert < 1 \right\} \subset \mathcal{R}_c
\end{equation}
It is worth mentioning that $\Omega_c$ is always a non-empty set. Indeed, for $\lambda_q=0$, we have $S(\theta,0,\Delta x)=1$ and consequently $\forall \Delta x \in \mathbb{R}_+$ there exists in $\mathbb{R}^2$ a neighborhood of $(\Delta x,\lambda_q)$ in which we have $\left\vert1-S(\theta,\lambda_q,\Delta x)\right\vert<1$. This means that we always have convergence for sufficiently small values of $\lambda_q$.
Lastly, we denote by $\mathcal{R}_m$ the region of stability of the modified equation : 
\begin{equation}
  \mathcal{R}_m  = \left\{ (\lambda_q,\Delta x) \in \mathbb{R}_+^2 : \left\vert\mathrm{e}^{\lambda_q \Delta x^qG(\theta,\lambda_q,\Delta x)}\right\vert \leq 1\right\} = \left\{ (\lambda_q,\Delta x) \in \mathbb{R}_+^2 : \mathrm{Re}(G(\theta,\lambda_q,\Delta x) \leq 0\right\}
\end{equation}
Provided these definitions, we can distinguish two cases : 
 \begin{enumerate}
     \item if $\mathcal{R}_s \subset \mathcal{R}_c$ then the stability of the modified equation provides reliable and complete information regarding the scheme stability. 
     \item If any subset of $\mathcal{R}_s$ lies outside of the convergence domain $\mathcal{R}_c$, this means that there is information on the stability limit of the scheme that is missed by the modified equations since its Fourier transform is non existing outside of $\mathcal{R}_c$. 
 \end{enumerate}
The following proposition is a direct consequence :

\begin{prpstn}
\label{pr:2}
For any scheme that writes as (\ref{eq:scheme}) and that is consistent with (\ref{eq:pde}) we have $\left(\mathcal{R}_m \cap \mathcal{R}_c \right)\subset \mathcal{R}_s$, that is if the modified equation is stable and its Fourier series is convergent then the scheme is also stable. 
\end{prpstn}
This proves that inside the convergence domain $\mathcal{R}_c$, the stability of the modified equation is a sufficient stability condition for the scheme. Furthermore, we shall add that if $\mathcal{R}_s \subset \mathcal{R}_c$ then this condition is also necessary. This result, literally, is not of practical interest as the stability of the full series $S(\theta,\lambda_q,\Delta x)$ is either nontrivial or impossible to obtain. But we will show that the above classification permits to justify whether a truncated version of the modified equation yields significant information on stability.


\subsection{Link between the stability of the scheme and the stability of a truncation}
Instead of the full series expansion, let us consider a truncated modified equation to an arbitrary order $N>P$: 
\begin{equation}
u_t = \sum_{p=1}^N \mu_p(\lambda_q,\Delta x) \frac{\partial^pu}{\partial x^p}, \quad \mu_N\neq 0
\label{eq:trunc}
\end{equation}
We recall that our main purpose is to know in which case does the stability of this truncation provide relevant information on the stability of the corresponding scheme. This differs from the approach of Warming and Hyett \cite{Fluid1974} in the sense that they showed how to reconstruct the exact Von Neumann amplification factor using a finite amount of coefficients $\mu_p$ without actually analyzing the stability of the truncated version. Under the Fourier transform, the previous equation writes : 
\begin{equation}
\frac{dv}{dt} =  \left( \sum_{p=1}^{N} \alpha_p(\lambda_q,\Delta x) \theta^p \right)v(k,t) =P_N(\theta,\lambda_q,\Delta x)v(k,t)
\label{eq:Fourier_mod_trunc}
\end{equation} 
Here, $P_N(\theta,\lambda_q,\Delta x)$ is a polynomial of $\theta$ of degree $N$ which is trivially a truncation of the series $G(\theta,\lambda_q,\Delta x)$. In the same manner as previously, the ordinary differential equation (\ref{eq:Fourier_mod_trunc}) yields :
\begin{equation}
v(k,t+\Delta t) = \mathrm{e}^{\Delta t P_N(\theta,\lambda_q,\Delta_x)}v(k,t) = S_N(\theta,\lambda_q,\Delta x)v(k,t)
\end{equation}
In contrast to the full series, the stability conditions of the truncation are obtainable in most cases through the analysis of the polynomial function $P_N(\theta,\lambda_q,\Delta x)$. Let $R_N(\theta,\lambda_q,\Delta x)$  be the rest of the series defined by : 
\begin{equation}
     R_N(\theta,\lambda_q,\Delta_x)=G(\theta,\lambda_q,\Delta_x)-P_N(\theta,\lambda_q,\Delta_x)=\sum_{p=N+1}^{\infty} \alpha_p(\lambda_q,\Delta x) \theta^p
\end{equation} 
In the convergence domain $\mathcal{R}_c$, the rest $R_N(\theta,\lambda_q,\Delta x)$ is bounded and we have : 
\begin{equation}
\lim\limits_{N \rightarrow +\infty} R_N(\theta,\lambda_q,\Delta x) \equiv 0 \quad \text{and} \ \lim\limits_{N \rightarrow +\infty} S_N(\theta,\lambda_q,\Delta x) \equiv S(\theta,\lambda_q,\Delta x)
\end{equation}
In this setting we can state the following result : 
\begin{prpstn}
Assume an initial condition satisfying $supp(v_0)\in [-M,M]$ and an arbitrary truncation order $N>P$, then for any $(\Delta x,\lambda_q)\in \mathcal{R}_c$, if the truncated modified equation is stable in the sense that there exists $C>0$ such that: 
\begin{equation}
\left\vert S_N(\theta,\lambda_q,\Delta x)\right\vert \leq 1 + C \Delta t
\end{equation}
then the scheme is also stable in the same sense.
\end{prpstn}
\begin{proof}
For $(\Delta x,\lambda_q)\in \mathcal{R}_c$ we have : 
\begin{equation*}
S(\theta,\lambda_q,\Delta x) = e^{\Delta t G(\theta,\lambda_q,\Delta x)} = S_N(\theta,\lambda_q,\Delta x)e^{\Delta t R_N(\theta,\lambda_q,\Delta x)}
\end{equation*}
Thus, since $v^n=S^n v_0$, we can write : 
\begin{align*}
\left\vert v^n \right\vert =\left\vert S^n v_0\right\vert &\leq \left\vert S_N \right\vert^n \left\vert e^{n\sum_{N+1}^\infty \alpha_p \theta^p}v_0\right\vert \\
&\leq (1+C\Delta t)^n \left\vert e^{nA\theta^{N+1}}v_0\right\vert \\
&\leq (1+C\Delta t)^n \left\vert e^{nA(k\Delta x)^{N+1}}v_0\right\vert \\
&\leq (1+C\Delta t)^n \left\vert e^{nA\Delta x^{N+1-q}\Delta tk^{N+1}/\lambda_q}v_0\right\vert  \\
&\leq e^{CT} \left\vert e^{AT\Delta x^{N+1-q}M^{N+1}/\lambda_q}v_0\right\vert \\
&\leq e^{CT}  e^{AT\Delta x^{N+1-q}M^{N+1}/\lambda_q}\left\vert v_0\right\vert 
\end{align*}
That is, $v^n$ is $L^2$-stable for initial conditions that are of compact support in the frequency domain, that is $k\in[-M,M]$.
\end{proof}

\section{Examples}
\subsection{Heat equation - centered finite differences}
\begin{figure}[H]
\includegraphics[]{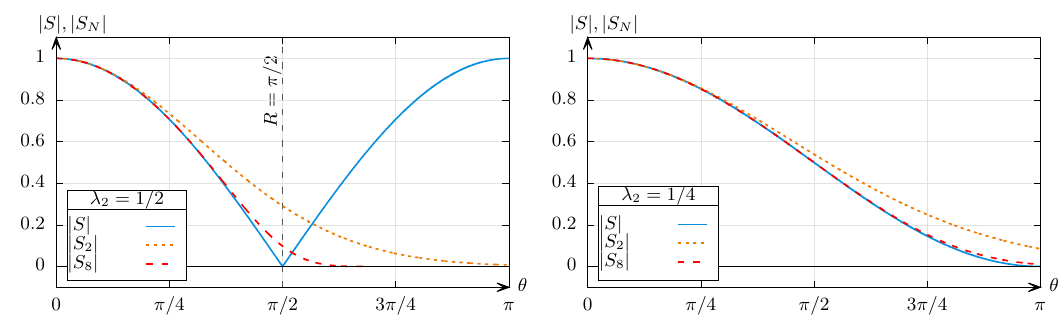}
\caption{Plot of the function $\left\vert S(\theta,\lambda_2,\Delta x)\right\vert$ along $\left\vert S_2(\theta,\lambda_2,\Delta x)\right\vert$ and $\left\vert S_8(\theta,\lambda_2,\Delta x)\right\vert$ for the values of $\lambda_2=1/2$ (left) and $\lambda_2=1/4$ (right). We can see that for $\lambda_2=1/2$, which lies outside of the convergence domain, the truncation curves stray away from the curve of $S$ starting from $\theta=R$. For $\lambda_2=1/4$, the truncations match well with $S$.} 
\label{fig:heat1}
\end{figure}

Consider the centered finite differences scheme for the heat equation. It can be cast into the form :   
\begin{equation}
    u_i^{n+1} = u_i^n + \alpha \lambda_2 \left(u_{i-1}^n - 2 u_i^n + u_{i-1}^n\right)
\end{equation}
We take for simplicity $\alpha=1$. The modified equations associated to this scheme up to $8th$ order for example is given by :
\begin{equation*}
   \frac{\partial u}{\partial t} = \frac{\partial^2u}{\partial x^2} + \frac{ \Delta x^2}{12}\left(1-6\lambda_2 \right) \frac{\partial^4u}{\partial x^4}+\frac{ \Delta x^4}{360}\bigg(1-30 \lambda_2\Big(1-4\lambda_2 \Big) \bigg)\frac{\partial^6u}{\partial x^6} +\frac{\Delta x^6}{20160} \Bigg(1-42 \lambda_2 \bigg(3 -40 \lambda_2 \Big( 1 -3 \lambda_2\Big)\bigg)\Bigg) \frac{\partial ^8u}{\partial x^8}+\ldots
\end{equation*}
Straightforward computations yield : 
\begin{equation}
S(\theta,\lambda_2,\Delta x) = 1 + 4 \lambda_2 \sin^2 (\theta/2)
\end{equation}
\begin{equation}
\mathcal{R}_s = \left\{ (\lambda_q,\Delta x) \in \mathbb{R}^2 :    \lambda_2 \leq \frac{1}{2} \right\} \quad ; \quad \Omega_c = \left\{ (\lambda_q,\Delta x) \in \mathbb{R}^2 :    \lambda_2 \leq \frac{1}{4} \right\} 
\end{equation}
\begin{equation}
\Delta t P_8(\theta,\lambda_2,\Delta x) = -\lambda_2 \left( \theta^2-\frac{1-6\lambda_2}{12} \theta^4 +\frac{1-30 \lambda_2\left(1-4\lambda_2 \right)}{360}\theta^6 -\frac{1-42 \lambda_2 \bigg(3 -40 \lambda_2 \Big( 1 -3 \lambda_2\Big)\bigg)}{20160}  \theta^8\right)
\end{equation} 
The domains of stability and convergence only depend on the parameter $\lambda_2$ independently of $\Delta x$. This permits to take an arbitrary value of $\Delta x=1$ and carry on the analysis, based on only $\lambda_2$. Furthermore since $\left\vert S(\theta,\lambda_2,\Delta x) \right\vert$ is an even function with respect to $\theta$ it suffices to look at $\theta \in [0,\pi]$.
Figure \ref{fig:heat1} shows a comparison between $\left\vert S(\theta,\lambda_2,\Delta x)\right\vert$ and $\left\vert S_N(\theta,\lambda_2,\Delta x)\right\vert$ for $N=2$ and $N=8$ in two cases. In the limit of stability $\lambda_2=1/2$, the scheme is stable but the series $G(\theta,\lambda_2,\Delta x)$ is not convergent for $\theta>R=\pi/2$ (See Appendix \ref{sec:append}). As displayed in the left-hand side of figure \ref{fig:heat1}, the curves begin aligned in the low frequencies, and then the truncation curves begin to diverge completely from the function $\left\vert S(\theta,\lambda_2,\Delta x)\right\vert$ once $\theta$ surpasses the threshold $R$. This is not the case for $\lambda_2=1/4$. For this value we can calculate the radius of convergence $R=\pi$ (See Appendix \ref{sec:append}). In fact, we can see on the right-hand side of the figure that for $\lambda_2$, the curves remain very close $\forall \theta \in [0,\pi]$. For $N=8$, the two curves almost overlap. 

\subsection{Transport equation - Upwind Euler}
The scheme writes : 
\begin{equation}
u_i^{n+1} = u_i^n - c \lambda_1 (u_i^n-u_{i-1}^n)
\end{equation}
We take $c=1$. In this case, the modified equation up to $4th$ order for example is given by : 
\begin{equation}
\frac{\partial u}{\partial t} +\frac{\partial u}{\partial x}= (1-\lambda_1)\left(\frac{\Delta x}{2}\frac{\partial^2 u}{\partial x^2} - \frac{\Delta x^2}{6}(1-2\lambda_1)\frac{\partial^3 u}{\partial x^3}+\frac{\Delta x^3}{24}(1-6\lambda_1(1-\lambda_1))\frac{\partial^4u}{\partial x^4}\right)+\ldots
\end{equation}
Straightforward computations yield : 
\begin{equation}
S(\theta,\lambda_1,\Delta x) = 1- \lambda_1\left(1-\mathrm{e}^{-\mathrm i \theta}\right)
\end{equation} 
\begin{equation}
\mathcal{R}_s = \left\{ (\lambda_q,\Delta x) \in \mathbb{R}^2 :    \lambda_1 \leq 1 \right\} \quad ; \quad \Omega_c = \left\{ (\lambda_q,\Delta x) \in \mathbb{R}^2 :    \lambda_1 \leq \frac{1}{2} \right\} 
\end{equation}
\begin{equation}
\Delta t P_4(\theta,\lambda_1,\Delta x)=-\mathrm{i}\lambda_1\theta+ \lambda_1(1-\lambda_1)\left(-\frac{1}{2}\theta^2 + \frac{1}{6}(1-2\lambda_1)\mathrm{i}\theta^3+\frac{1}{24}(1-6\lambda_1(1-\lambda_1))\theta^4\right)+\ldots
\end{equation}

In contrast to the previous example for the heat equation, we could not compute explicitly the radius of convergence for values of interest of $\lambda_1$. Nevertheless, it is possible to extend $\Omega_c$ to cover all the values of $\lambda_1$ that are in $\mathcal{R}_s$ (See appendix \ref{app:B}). Hence, as shown in figure \ref{fig:trans1}, for values of $\lambda_1\leq 1$, the amplification factor of the truncation to $6th$ order $\left\vert S_6(\theta,\lambda_1,\Delta x) \right\vert$ seems to be a very good approximation of the scheme amplification factor $\left\vert S(\theta,\lambda_1,\Delta x) \right\vert$, even for values of $\lambda_1$ that are in the vicinity of the stability threshold. 
\begin{figure}[H]
\includegraphics[]{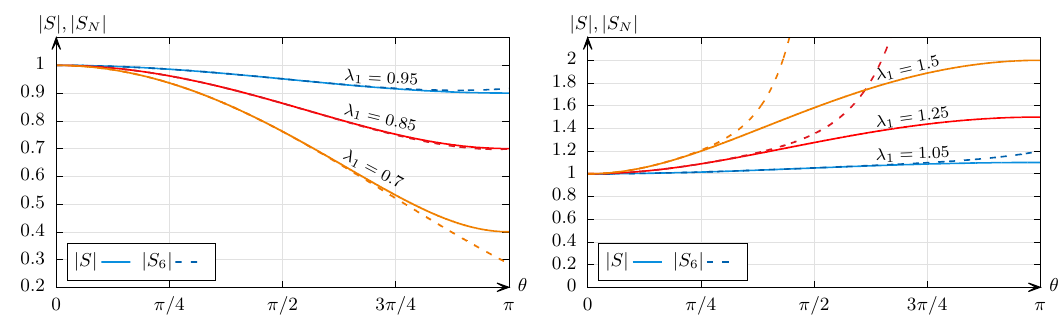}
\caption{Plot of the function $\left\vert S(\theta,\lambda_1,\Delta x)\right\vert$ along $\left\vert S_6(\theta,\lambda_2,\Delta x)\right\vert$ for different values of $\lambda_1$ in the stable region (left) and in the unstable region (right). The continuous lines stand for the amplification factor $\left\vert S(\theta,\lambda_1,\Delta x)\right\vert$ and the dashed lines represent the truncated modified equation amplification factor $\left\vert S_6(\theta,\lambda_1,\Delta x)\right\vert$.  We can see that for $\lambda_1\leq1$ the truncation curves match well with the exact amplification factor even in the boundaries of stability. Inversely, it seems according to the right-hand graphic that $\lambda_1=1$ marks the threshold of convergence. } 
\label{fig:trans1}
\end{figure}
\section*{Conclusion and perspectives}
We could explain throughout this work that one of the main reasons behind the failures of the modified equations technique is non other than series divergence. Although the analysis only provides sufficient conditions in general, it lifts some well-known ambiguities and provides some assumptions required by the technique to be of a more justified practical use. We show that in these settings, the stability of a truncation gives already reasonable approximations to stability conditions of the scheme. An extension of the obtained results to the case of systems is underway. It would be also interesting to extend this approach to cover a larger class of explicit and also implicit schemes.    
\appendix
\section{Computing some radii of convergence}
\label{sec:append}
In each of the following cases, we would like to check the convergence radius of the series : 
\begin{equation}
\ln(S(\theta,\lambda_q,\Delta x) = \sum_{p=0}^\infty \frac{\alpha_p}{\lambda_q \Delta x^q} \theta^p = \sum_{p=0}^\infty a_p \theta^p
\end{equation}
For that, we use the root convergence test, that is the radius of convergence $R$ is given by :
\begin{equation}
\lim\limits_{p \rightarrow +\infty} \left\vert a_p\right \vert ^{1/p}=\frac{1}{R}
\end{equation}

\subsection*{Centered scheme for the heat equation : $\lambda_2=\frac{1}{2}$}
We have :
\begin{equation}
S(\theta,\lambda_2,\Delta x) = 1-2\sin(\theta/2)^2
\end{equation}
In this case, it is possible to obtain an explicit expression of the $pth$ term of the sequence $a_p$ which is given by : 
\begin{equation}
\left\{\begin{array}{ll}
a_0=0 \\
\displaystyle a_{2p} = -\frac{(-4)^pE_{2p-1}(0)}{2(2p)!}  &\displaystyle \forall p\geq 1 \\
a_{2p+1} = 0 \quad &\forall p\geq 0
\end{array} \right.
\end{equation}
where $E_p(x)$ denotes the $pth$ Euler polynomial function defined by : 
\begin{equation}
\frac{2 e^{xt}}{e^t+1} = \sum_{n=0}^\infty E_n(x) \frac{t^n}{n!}
\end{equation}
Since the series has only even order coefficients, it is equivalent and more convenient to check the convergence of the series : 
\begin{equation}
\sum_{p=0}^\infty b_p \phi^p
\label{eq:phi}
\end{equation}
where $\phi=\theta^2$ and $b_p=a_{2p}$. We proceed as follows in order to check the radius of convergence. We use Bernoulli's number \cite{euler}:
\begin{equation}
B_{2p} = \frac{-p E_{2p-1}(0)}{(2^{2p}-1)}
\end{equation} 
and plug this ansatz into $b_p$ to obtain : 
\begin{align*}
\left\vert b_p \right\vert= \left\vert\frac{(4)^pE_{2p-1}(0)}{2(2p)!}\right\vert = \left\vert\frac{(4)^p(2^{2p}-1)}{2p(2p)!}B_{2p} \right\vert\underset{p\rightarrow\infty}\sim\left\vert\frac{2^{2p}}{p(\pi)^{2p}}\right\vert
\end{align*}
Hence we can write : 
\begin{equation}
\left\vert b_p\right \vert \underset{p\rightarrow\infty}\sim \frac{1}{p}\left(\frac{2}{\pi}\right)^{2p}
\end{equation}
and so: 
\begin{equation}
\lim\limits_{p \rightarrow +\infty} \left\vert b_p\right \vert ^{1/p}=\left(\frac{2}{\pi}\right)^2
\end{equation}
This means that the series (\ref{eq:phi}) has a radius $R_{\phi}=\pi^2/4$ which implies that the radius of convergence of the series $G(\theta,1/2,\Delta x)$ is $R=\pi/2$.
\subsection*{Centered scheme for the heat equation : $\lambda_2=\frac{1}{4}$}
We have :
\begin{equation}
S(\theta,\lambda_2,\Delta x) = 1-\sin(\theta/2)^2
\end{equation}
The expression of the $pth$ term $a_p$ is given by : 
\begin{equation}
\left\{\begin{array}{ll}
a_0=0 \\
\displaystyle a_{2p} = -\frac{(-1)^pE_{2p-1}(0)}{(2p)!}  &\displaystyle \forall p\geq 1 \\
a_{2p+1} = 0 \quad &\forall p\geq 0
\end{array} \right.
\end{equation}
Using the same previous notations as in (\ref{eq:phi}), we have : 
\begin{equation}
\left\vert b_p \right\vert = \left\vert\frac{(E_{2p-1}(0)}{(2p)!} \right\vert= \left\vert\frac{(2^{2p}-1)}{p(2p)!}B_{2p}\right\vert\underset{p\rightarrow\infty}\sim \left\vert\frac{2}{p\pi^{2p}}\right\vert
\end{equation}
which yields :
\begin{equation}
\lim\limits_{p \rightarrow +\infty} \left\vert b_p\right \vert ^{1/p}=\frac{1}{\pi^2}.
\end{equation}
This gives the radius of convergence of the series $G(\theta,1/4,\Delta x)$ is $R=\pi$.
\section{Proof of convergence for $\lambda_1 \leq 1$ of Upwind Euler for transport equation}
\label{app:B}
Since $\Omega_c = \left\{ (\lambda_q,\Delta x) \in \mathbb{R}^2 :    \lambda_1 \leq \frac{1}{2} \right\}$, then the series $G(\theta,\lambda_1,\Delta x)$ is convergent $\forall \lambda_1 \leq 1/2$ and we have : 
\begin{equation*}
e^{\Delta t G(\theta,\lambda_1,\Delta x)} = S(\theta,\lambda_1,\Delta x)
\end{equation*} 
Since stability only depends on the modulus of $S(\theta,\lambda_1,\Delta x)$ it is sufficient for our analysis to consider the equality : 
\begin{equation}
e^{\Delta x \lambda_1 \mathrm{Re}(G(\theta,\lambda_1,\Delta x))} = \left\vert S(\theta,\lambda_1,\Delta x) \right\vert
\label{eq:real}
\end{equation} 
Next, we show the following symmetry : 
\begin{equation*}
\left\vert S(\theta,1/2-\lambda_1,\Delta x) \right\vert = \left\vert S(\theta,1/2+\lambda_1,\Delta x) \right\vert
\end{equation*}
Indeed we have : 
\begin{align*}
S(\theta,1/2-\lambda_1,\Delta x) = 1- (1/2-\lambda_1)\left(1-\mathrm{e}^{-\mathrm i \theta}\right) &= (1/2+\lambda_1) + (1/2-\lambda_1)e^{-\mathrm i\theta} \\
S(\theta,1/2+\lambda_1,\Delta x) = 1- (1/2+\lambda_1)\left(1-\mathrm{e}^{-\mathrm i \theta}\right) &= (1/2-\lambda_1) + (1/2+\lambda_1)e^{-\mathrm i\theta}  \\
&= \left((1/2+\lambda_1) + (1/2-\lambda_1)e^{\mathrm i\theta}\right)e^{-\mathrm i\theta} \\
&= e^{-\mathrm i\theta}\overline{S(\theta,1/2-\lambda_1,\Delta x)}
\end{align*}
where the bar denotes the complex conjugate. Hence $\left\vert S(\theta,1/2-\lambda_1,\Delta x) \right\vert = \left\vert S(\theta,1/2+\lambda_1,\Delta x) \right\vert$. This implies through equality (\ref{eq:real}) that : 
\begin{equation}
\Delta x (1/2-\lambda_1) \mathrm{Re}(G(\theta,1/2-\lambda_1,\Delta x)) = \Delta x (1/2+\lambda_1) \mathrm{Re}(G(\theta,1/2+\lambda_1,\Delta x))
\end{equation}
which also implies in terms of coefficients: 
\begin{equation}
(1/2-\lambda_1)\alpha_{2p}(1/2-\lambda_1,\Delta x) = (1/2+\lambda_1)\alpha_{2p}(1+\lambda_1,\Delta x) \quad \forall p\geq 1
\end{equation}
Therefore, since the series is convergent for $0\leq\lambda_1 \leq 1/2$, it is also convergent for $1/2\leq \lambda_1 \leq 1$. 
\bibliography{biblio}
\bibliographystyle{unsrt}  
\end{document}